\documentclass[twoside,10pt]{article}

\usepackage{amssymb,amsmath,latexsym,epsfig,euscript,theorem}

\font\tengoth=eufm10 \font\sevengoth=eufm7 \font\fivegoth=eufm5
\newfam\gothfam  \textfont\gothfam=\tengoth
\scriptfont\gothfam=\sevengoth \scriptscriptfont\gothfam=\fivegoth

\def\rank {{\rm \mbox{rank}}}

\def\Nor{{\rm \mbox{Nor\,}}}
\def\Ver{{\rm \mbox{Ver\,}}}

\def \R{\mathbb{R}}
\newcommand{\C}        {\mathbb{C}}

\def \P{\mathbb{P}}

\def \sp {\mathbb{S}}
\def \m {\overline m}

\def \P{\mathbb{P}}

\newtheorem{defn}{Definition}[section]

\newtheorem{prop}[defn]{Proposition}

\newtheorem{thm}[defn]{Theorem}

{\theorembodyfont{\rmfamily} \newtheorem{rem}[defn]{Remark}}

\newtheorem{conj}[defn]{Conjecture}

\def\nna {\nabla^{\perp}}

\def\Ka {K\"ahler }

\def\norp {N_p(M)}

\def\rank {{\rm rank \; }}
\def\span {{\rm span \; }}

\def\Ph {\Phi^{\perp}}

\def\p {\mathfrak{p}}
\def\m {\mathfrak{m}}
\def\k {\mathfrak{k}}
\def\g {\mathfrak{g}}

\begin{document}

\title{Parallel submanifolds of complex projective space and
their normal holonomy\footnote{Research partially supported by GNSAGA (INdAM) and MIUR of Italy}}

\author{ Sergio Console and Antonio J. Di Scala}

\date{}

\maketitle

\begin{abstract}
The object of this article is to compute the holonomy group of the
normal connection of complex parallel submanifolds of the complex
projective space. We also give a new proof of the classification
of complex parallel submanifolds by using a normal holonomy
approach. Indeed, we explain how these submanifolds can be
regarded as the unique complex orbits of the (projectivized)
isotropy representation of an irreducible Hermitian symmetric
space. Moreover, we show how these important submanifolds are related
to other areas of mathematics and theoretical physics. Finally, we
state a conjecture about the normal holonomy group of a complete
and full complex submanifold of the complex projective space.
\end{abstract}


\vspace{2mm}


\date{}

\maketitle

\section{Introduction.}

The study of the normal holonomy group, started by C. Olmos in
\cite{Olmos} (see also \cite{BCO} for more details and
applications), turned out to be a powerful tool
for the study of submanifolds with simple geometric invariants,
e.g. homogeneous submanifolds, isoparametric submanifolds and
their generalizations \cite{Olmos4,DiOl}. In particular, normal
holonomy methods were used by C. Olmos to give simpler and
geometric proofs of Berger-Simons' theorems on holonomy
\cite{Olmos2,Olmos3}.

The restricted normal holonomy group $\Ph_p$ of a complex
submanifold $M \subset \C P^N$ at a point $p \in M$ acts on the
normal space $\norp$. Under suitable and very general conditions
(see \cite{AlDi} and Section~\ref{prelim}) this action agrees with
the isotropy action of an irreducible Hermitian symmetric space;
i.e., the pair $(\Ph_p,\norp)$ is given by $(K, T_{[K]} G/K) = (K,
{\p})$, (with ${\g} = \k \oplus \p$ a Cartan decomposition) where
$G/K$ is an irreducible Hermitian symmetric space and the action
coincides with the isotropy representation of  $G/K$ on $\p$.

Symmetric submanifolds play a distinguished role among
submanifolds with simple geometric invariants. They are analogous
to symmetric spaces for submanifold theory. Indeed, they always
come equipped with a symmetry at each point, namely the geodesic
reflection in the corresponding normal submanifold. This implies
that the second fundamental form is parallel. Actually, for a
complex submanifold of $\C P^N$ being extrinsically symmetric is
the same as having parallel second fundamental form (see
\cite[Proposition 9.3.1, page 256]{BCO} for more details). For
real space forms, symmetric submanifolds turn out to be orbits of
the isotropy representation, the so called $R$-spaces. This is a
classical result of D. Ferus \cite{Fer} (see also \cite[Chapter 3,
Section 7]{BCO}). Notice that these submanifolds are also related
to the theory of Jordan algebras \cite[p. 239]{Ber}.

In algebraic geometry, complex symmetric submanifolds of $\C P^N$
are called \it characteristic projective subvarieties\rm. They
are related with (a part of) the celebrated Borel-Weil theorem
since they are the unique complex orbit of the action of a compact
Lie group (see \cite{BoWe}, \cite[p. 166]{GuSt} and \cite[p. 103, Remark]{Mok}).
Furthermore they are the main ingredient of a polarization-type
argument used by N. Mok to prove his well-known rigidity theorems
for higher rank Hermitian symmetric spaces \cite[p. 111, Prop.
3]{Mok}. In the theory of Jordan algebras, complex symmetric
submanifolds of $\C P^N$ are described by means of \it minimal
tripotents, \rm see \cite[p. 579]{Kaup}. These submanifolds are
also important in physics and chemistry. Namely, they are related
to the so-called approximation of Hartree-Fock  (e.g. Slater
determinants, for details see \cite[p. 165]{GuSt}).

Complex symmetric submanifolds of $\C P^N$ were classified by
Nakagawa-Takagi \cite{NaTa}. Besides the standard techniques from
representation theory (cf. \cite{Take}), their
classification depends upon the work of Calabi-Vesentini and Borel
\cite{CaVe, Borel}. Namely, they managed to give a link between
the norm of the covariant derivatives of the second fundamental
form of the canonical embeddings of Hermitian symmetric spaces
and the eigenvalues of the curvature operator introduced by Nakano
\cite{Na}.

The first aim of this paper is to compute the pairs
$(\Ph_p,\norp)$ when $M$ is a complex parallel submanifold of the
projective space $\C P^N$. We are going to collect our
results in the third column of Table 1  below. Besides the
classification of the possible pairs $(\Ph_p,\norp)$ given in
\cite{AlDi}, our method is based on Proposition \ref{slice} and
the classification of complex parallel submanifolds.

\begin{table} [h!]\ {\footnotesize
\vskip .5cm {}
\begin{tabular}{ | c | c | c | c | }\hline
\rule[-20pt]{0pt}{33pt}%
\parbox{3cm}{\center{Hermitian symmetric space $G/K$}}  & \parbox{2cm}{\center{$M$ as complex $K$-orbit}}
& \parbox{1.5cm}{\center {Normal holonomy}} &
\parbox{1.5cm}{\center {Remarks}}\\ \hline\hline
\rule [-14pt]{0pt}{33pt}%
$\dfrac{E_7}{T^1 \cdot E_6}$ & $\dfrac{E_6}{T^1 \cdot Spin_{10}}$
& $\dfrac{SO(12)}{T^1 \cdot SO(10)}$ & \\ \hline
 \rule [-14pt]{0pt}{33pt}%
$\dfrac{E_6}{T^1 \cdot Spin_{10}}$ & $\dfrac{SO(10)}{U(5)}$  &
$\dfrac{U(6)}{U(5)}$ &  \\ \hline
 \rule [-14pt]{0pt}{33pt}%
 $\dfrac{Sp(n+1)}{U(n+1)}$ & $\C P^n$  & $\dfrac{Sp(n)}{U(n)}$& Veronese
  \\ \hline
 \rule [-14pt]{0pt}{33pt}%
 $Gr_2^+(\R^{n+2}):=\dfrac{SO(n+2)}{T^1 \cdot SO(n)}$ & $Gr_2^+(\R ^{n})$  & $\dfrac{U(2)}{U(1)}$ & Quadrics
  \\ \hline
 \rule [-14pt]{0pt}{33pt}%
 $\dfrac{SO(2n)}{U(n)}$ & $Gr_2(\C ^n)$  & $\dfrac{SO(2(n-2))}{U(n-2)}$& Pl\"ucker \\ \hline
 \rule [-14pt]{0pt}{33pt}%
$Gr_a(\C^{a+b}):=\dfrac{SU(a+b)}{S(U(a)\times U(b))}$ & $\C
P^{a-1} \times \C P^{b-1}$  & $\dfrac{SU(a+b - 2)}{S(U(a-1)\times
U(b-1))}$ & Segre \\ \hline
\end{tabular}
}

\caption{\footnotesize Symmetric complex submanifolds $M \subset \P
(T_{[K]}G/K)$. The space in the third column is the Hermitian symmetric space whose
 isotropy representation gives the normal holonomy action.}
\end{table}

The second goal of this paper is to obtain the classification of
complex parallel submanifolds of the complex projective space
without making use of the work of Calabi-Vesentini and Borel
\cite{CaVe, Borel}, as in the classical work of Nakagawa and
Takagi \cite{NaTa}. Indeed, such a classification is contained in
the second column of Table 1. It turns out that complex parallel
submanifolds of $\C P^n$  consist of the first canonical
embeddings of  rank two Hermitian symmetric spaces, the Veronese
and Segre embeddings. Our approach (presented in
Section~\ref{classification})
is based
on holonomy techniques and the knowledge of the codimension
of the canonical embeddings only.
We will also explain how complex parallel submanifolds of $\C P^N$
can be regarded as the unique complex orbits of the (projectivized)
isotropy representation of an irreducible Hermitian symmetric
space. Notice the analogy with the above mentioned result of Ferus in the real
setting.

Our main result is encoded in Table 1 and can be read as follows.\\
\textit{Let $G/K$ be an irreducible Hermitian symmetric space from
the complete list  in the first column of Table 1. Then, the
second column contains the unique complex orbit of the $K$-action
on the projective space $\P (T_{[K]}G/K)$. This $K$-orbit has
parallel second fundamental form and all complex parallel
submanifolds arise in this way. Moreover, the third column of the
table contains the Hermitian symmetric space whose isotropy
representation gives the normal holonomy action $(\Phi_p^{\perp},
N_p(M))$.}

In particular we prove the following

\begin{thm} \label{main} Let $M \subset \C P^N$ be a full (connected)
complex submanifold with parallel second fundamental form. Then
$M$ is an open subset of the unique complex orbit of the
(projectivized) isotropy representation of an irreducible
Hermitian symmetric space $G/K$. Moreover, their normal holonomy
actions $(\Phi_p^{\perp}, N_p(M))$ agree with the isotropy
representations of the Hermitian symmetric spaces listed in the
third column of Table 1. \end{thm}

For higher canonical embeddings of Hermitian symmetric spaces we
get

\begin{thm} \label{alto} Let $f_d: G/K \hookrightarrow \C P^{N_d}$ be the
$d-$th canonical embedding of an irreducible Hermitian symmetric
space. If $d > 2$ then the normal holonomy group is the full
unitary group of the normal space.
\end{thm}

Motivated by the above theorem we propose the following extrinsic
analog of Berger's theorem as a conjecture

\begin{conj} \label{ProBer} Let $M \hookrightarrow \C P^N$ be a complete (connected) and full (i.e. not
contained in a proper hyperplane) complex submanifold. If the
normal holonomy group is not the full unitary group, then $M$ has
parallel second fundamental form.
\end{conj}

Notice that if the above conjecture is true, then the realization
problem of normal holonomy group of complex submanifolds of $\C
P^N$ is solved. Namely, up to the isotropy representation of the
exceptional $\dfrac{E_7}{T^1 \cdot E_6}$ any other isotropy
representation of an irreducible Hermitian symmetric space can be
obtained as a normal holonomy action. Recall that the realization problem of the normal
holonomy group of
submanifolds of the sphere  was solved in \cite{HeOl}, up to eleven exceptions.
Finally, Conjecture \ref{ProBer} can be regarded as the
complex version of the conjecture posed in \cite{Olmos5}. Namely,
\it an irreducible and full homogeneous submanifold of the sphere,
different from a curve, whose normal holonomy group is not
transitive, must be an orbit of an s-representation. \rm

\section{Preliminaries.}\label{prelim}

Throughout this paper by \it complex submanifold $M$ of $\C P^N$
\rm we mean a holomorphic and isometric embedding $M
\hookrightarrow \C P^N$, where $\C P^N$ carries the standard
Fubini-Study K\"ahler form of constant holomorphic curvature $1$.
We will always assume that $M$ is connected.

We refer to  \cite{BCO} for the definitions of the normal bundle
$N(M) \to M$, its normal connection and its holonomy group $\Ph_p$
at a point $p \in M$ (the so called \it normal holonomy group\rm).

Recall that the shape operator of a complex submanifold anticommutes with
the complex structure $J$, i.e., $A_\xi J=J A_\xi$, for any normal vector $\xi$.
Moreover, the equation of Gauss yields the following expression for the holomorphic sectional curvature of $M$
$$
\frac 12 \left (\langle X, Y \rangle^2 + \langle X, JY \rangle^2+ \Vert X \Vert^2\Vert Y \Vert^2 \right )
 = \langle R_{X, JX} JY, Y \rangle+ 2  \Vert \alpha (X, Y) \Vert^2\, . \leqno {(*)}
$$

We say that $M \subset \C P^N$ is \textit{full} if it is not contained in
a hyperplane of $\C P^N$. Recall that the first normal space
$N^1(M)$ is defined by $N^1(M) := \span \{\alpha(X,Y)\}$.

The following result gives a sufficient condition for a submanifold of $\C P^N$ not to be full.

\begin{thm} {\rm \cite{Ce}, \cite{ChO1}} \label{codi}
Let $M$ be a \Ka submanifold of $ \C P^N$. If there exists a
complex $\nna$-parallel subbundle $V \neq 0$ of the normal bundle
$N(M)$ such that $V \perp \, N^1(M)$, then $M$ is non-full.
\end{thm}

Calabi rigidity theorem of complex submanifolds $ M \hookrightarrow \C
P^N$ \cite{Cal} implies that isometric and holomorphic immersions
are equivariant (see Subsection \ref{canonical} ) Namely, any
intrinsic isometry can be extended to the ambient space, i.e., to
the projective space $\C P^N$. For a detailed explanation see
\cite[p. 655, Theorem 4.3]{NaTa} or \cite{Take}.

Finally, let us recall the following well-known proposition (see e.g. \cite{NaTa}).

\begin{prop} \label{complete} Let $M$ be a K\"ahler manifold not necessarily complete.
Let $f: M \rightarrow \C P^N$ be a holomorphic and isometric
immersion with parallel second fundamental form. Then $M$ is a
locally Hermitian symmetric space of compact type. Moreover, there
exists a complete Hermitian symmetric space of compact type
$\widetilde{M}$ and a holomorphic and isometric embedding
$\widetilde{M} \stackrel{\widetilde{f}}\hookrightarrow \C P^N$with
parallel second fundamental form such that $f = \widetilde{f}
\circ i $, where $i:M \rightarrow \widetilde{M}$ is the canonical
inclusion.
\end{prop}

\subsection{Normal holonomy.}
The link between isotropy and holonomy is well-known for
Riemannian symmetric spaces. For submanifolds with
parallel second fundamental form there is a similar relationship
between isotropy and normal holonomy group.

As we quoted in the introduction, it was proved in \cite{AlDi}
that, under suitable and very general conditions, the normal holonomy group acts
on the normal space as the isotropy representation of an irreducible Hermitian symmetric
space. All known
examples (complete submanifolds,
K\"ahler-Einstein submanifolds and manifolds with zero index of
relative nullity) satisfy these conditions. Moreover, it is not known whether there exists a
full complex submanifold whose normal holonomy group acts
in reducible way. For completeness let us state here the following special
case of the main result in \cite{AlDi}.

\begin{thm} {\rm \cite{AlDi}} \label{holno} Let $M \subset \C P^N$ be
a complete and full submanifold of $\C P^N$. Let $\Ph_p$ be the
normal holonomy group at $p \in M$. Then, there exists an
irreducible Hermitian symmetric space $H/S$ such that $\Ph_p = S$.
Indeed, $\dim_{\C}(N_p(M)) = \dim_{\C}(H/S)$ and the
$\Ph_p$-action on $N_p(M)$ agrees with the isotropy representation
of $S$ on $T_{[S]}(H/S)$.
\end{thm}

The following proposition gives a nice application of the above
theorem.

\begin{prop} \label{slice} Let $M=G/K \hookrightarrow \C P^N$ be a full parallel
submanifold where $M=G/K$ is a Hermitian symmetric space. Then the
normal holonomy group $\Ph_p$ is homomorphic with $K$ i.e. $\Ph_p
= K/I$, where $I$ is normal in $K$. \end{prop}

\it Proof. \rm Let us denote by  $K^{\perp}$ the image of the restriction of the isotropy
representation of $K$ to the normal space $N_p(M)$,  the so
called \it slice representation \rm. Thus, we are going to show
that $K^{\perp} = \Ph_p$. First of all, notice that,  since  isometries preserve parallel transport, $K^{\perp}
\subset \Nor(\Ph_p)$, where $\Nor(\Ph_p) \subset U(N_p(M))$ is the
normalizer of $\Ph_p$ in the full unitary group $U(N_p(M))$.
 By the
above Theorem \ref{holno} $\Ph_p$ is isomorphic to the isotropy
$S$ of an irreducible Hermitian symmetric space. Thus, $\Nor(\Ph_p)
= \Ph_p$ and we get the inclusion $ K^{\perp} \subset \Ph_p$.

The proof that $\Ph_p \subset K^{\perp}$ is similar to the one
given in \cite[p.7, Theorem 2]{Esch}. Namely, any transvection of
$G/K$, when extended to the ambient projective space, gives
the parallel transport with respect to the normal connection. Then
we can approximate any closed curve by a closed geodesic polygon.
So we get a composition of isometries which belong to $K$ and, by
construction, to $\Ph_p$.\\ By taking limits and using the
compactness of the involved groups we get the desired inclusion
$ \Ph_p \subset K^{\perp}$. {}\hfill $\Box$

\vspace{.5cm}

As a sum up of the above propositions, we have the following theorem.

\begin{thm} \label{norhol} Let $M=G/K \hookrightarrow \C P^N$ be a
Hermitian symmetric space embedded into $\C P^N$ with parallel
second fundamental form. Assume also that the embedding is full.
Then, there exists an irreducible Hermitian symmetric space $H/S$
such that $\Ph_p = S = K/I$ where $I \subset K$ is a normal
subgroup, $\dim_{\C}(N_p(M)) = \dim_{\C}(H/S)$ and $\Ph_p$ acts on
$N_p(M)$ as the isotropy representation of $S$ on $T_{[S]}(H/S)$.
\end{thm}

\subsection{Parallel products.}

The following result of Nakagawa and Takagi allows us to restrict to parallel
embeddings of irreducible Hermitian symmetric spaces. Here we present an alternative proof.

\begin{thm} {\rm \cite[p. 664, Lemma 7.1]{NaTa}} Let $M_i$ be an
$n_i$-dimensional K\"ahler manifold $(i=1,2)$. If the K\"ahler
manifold $M_1 \times M_2$ admits a K\"ahler immersion into $\C
P^{n_1 + n_2 + p}$ with parallel second fundamental form, then
$M^i$ is locally $\C P^{n_i}$  $(i=1,2)$.
\end{thm}

\it Proof. \rm We can assume that the immersion is full otherwise
by using Theorem \ref{codi} we can reduce the codimension. So let
us introduce the subbundles $\alpha_{1,1}, \, \alpha_{2,2}$ and
$\alpha_{1,2}$ of the normal bundle $N(M)$ as follows: \[
\alpha_{i,j} := \{ \alpha(TM_i,TM_j) \}, \] where $\alpha$ is the
second fundamental form of $M$. Since the immersion $M_1 \times
M_2 \subset \C P^{n_1 + n_2 + p}$ has parallel second fundamental
form, we get  $N(M) = \alpha_{1,1} + \alpha_{1,2} + \alpha_{2,2}$.
Indeed, this is a consequence of Theorem \ref{codi} since the
first normal space is parallel. Moreover, a simple application of
the equation $(*)$ implies that the above sum is orthogonal, i.e.
$N(M) = \alpha_{1,1} \oplus \alpha_{1,2} \oplus \alpha_{2,2}$.
Observe that any $\alpha_{i,j}$ is a parallel subbundle with
respect to the normal connection. Thus, since $\Ph_p$ acts
irreducibly on $N_p(M)$ by Theorem \ref{holno}, two of the three
subbundles $\alpha_{i,j}$ must be trivial. It is not difficult to
check that $\alpha_{1,2}$ cannot be trivial (see \cite[p. 202,
Thm. 17]{AlDi} for details). Indeed, equation $(*)$ with $X \in
TM^1$ and $Y \in TM^2$ yields a contradiction. Using again the
equation $(*)$, we get that the curvature tensor of each factor
agrees with the curvature tensor of the ambient space $\C P^{n_1 +
n_2 + p}$ and we are done. \hfill $\Box$

\subsection{Canonical embeddings.}\label {canonical}

Let us assume now that $M=G/K$ is an irreducible Hermitian
symmetric space and let $f : M \hookrightarrow \C P^N$ be a full
holomorphic and isometric embedding. From Calabi rigidity theorem
it follows that the embedding $f : M \hookrightarrow \C P^N$ is
$G$-equivariant (see \cite[p. 655, Theorem 4.3]{NaTa} or
\cite{Take}). Then, from Elie Cartan's work, such embeddings $G/K
\hookrightarrow \C P^N$ are well-known and are called
\textit{canonical embeddings}.

They can be constructed by means of the representation theory of the
simple group $G$ through the so-called Borel-Weil construction
(see \cite{BoWe, Take}), which holds, more generally, for homogeneous K\"ahler
manifolds  and can be summarized as follows.

Let $d$ be a positive integer and $\rho : G^\C \to
\mathfrak{gl} (\C^{N_d+1})$ the irreducible representation of the
complexification $G^\C$ of $G$ with highest weight $d\Lambda_j$,
where $\Lambda_j$ is the fundamental weight corresponding to the
simple root $\alpha_j$. Let $p$ be a highest weight vector
corresponding to  $d\Lambda_j$. Then the action of $G^\C$ on
$\C^{N_d+1}$ induces a unitary representation of $G$ whose orbit
of the highest weight vector $p$ in  $\C P^{N_d}$ yields  a full
holomorphic embedding $f_d$ of $M=G/K$ into $\C P^{N_d}$,
called the {\it d-th canonical embedding} of $M$ into a complex
projective space.

The submanifold $M$ of $\C P^{N_d}$ is the
unique complex orbit of the action of  $G$ on $\C P^{N_d}$ (or
equivalently, the unique compact orbit of the $G^\C$-action). The
dimension $N_d$ can be calculated explicitly by means of the Weyl's
dimension formula. The induced metric on $M \subset \C P^{N_d}$ is
K\"ahler-Einstein.


\section{Normal holonomy of parallel submanifolds.}\label{normal_hol}

In this section we are going to prove the last sentence of
Theorem~\ref{main} and later Theorem~\ref{alto}. The proof of the
first part of Theorem~\ref{main} will be given in
Section~\ref{classification}.

\it Proof of the last part of
Theorem~\ref{main}. \rm
We compute the third column of Table 1 using Theorem \ref{norhol}.

\medskip

Let us start with the first line, explaining our method. Namely,
let us focus on the $1-$st canonical embedding $\dfrac{E_6}{T^1
\cdot Spin_{10}} \hookrightarrow \C P^{26}$. By Theorem
\ref{norhol}, the normal holonomy group $\Ph_p $ is a quotient of $
T^1 \cdot Spin_{10}$ whose action on the $10$-dimensional normal
space agrees with the isotropy representation of a $10$-dimensional
irreducible Hermitian space. Looking for  such an  irreducible
Hermitian symmetric space (see the first column of Table 1) we see
that the only possibility is $\dfrac{SO(12)}{T^1 \cdot SO(10)}$.
Thus, $\Ph_p $ acts on $N_p(M)$ as the isotropy representation of
$\dfrac{SO(12)}{T^1
\cdot SO(10)}$.\\
The computation for the second line is similar. Indeed, by
Theorem \ref{norhol}, the normal holonomy group $\Ph_p $ of the
embedding $\dfrac{SO(10)}{U(5)} \hookrightarrow \C P^{15}$ is a
quotient of $U(5)$ whose action on the $5$-dimensional normal
space agrees with a isotropy representation of a $5$-dimensional
irreducible Hermitian space. Looking for such an irreducible
Hermitian symmetric space,  the only possibility is
$\dfrac{U(6)}{U(5)}$. Thus,  $\Ph_p $ acts on $N_p(M)$ as the
isotropy representation of $\dfrac{U(6)}{U(5)}$, i.e. as the standard representation of $U(5)$ on $\C^5$.\\
For the next four classical cases (i.e. the Veronese, the
Quadrics, etc) a similar analysis, based on Theorem \ref{norhol},
can be done in order to compute the third column. This completes
the proof of the last sentence in Theorem \ref{main}. \hfill
$\Box$

\vspace{.5cm}

\it Proof of Theorem~\ref{alto}. \rm It is enough to prove that
the first normal space does not agree with the full normal space.
We will actually show
$$
\dim_{\C} N_p > \frac{m(m+1)}2 \, , \qquad m=\dim_{\C} (G/K)\, ,\leqno(**)
$$
which implies the above assertion, since the dimension of the
first normal space is smaller or equal than $\dfrac{m(m+1)}2$.
Indeed, if the first normal space has a complement, any vector in
this complement belongs to the nullity of the adapted normal
curvature tensor $\tilde {\cal R}$ (see \cite{AlDi, Olmos}). This
is absurd, since the curvature tensor of a Hermitian symmetric
space has no nullity, so that the normal holonomy must be
transitive on the unit sphere of the normal space.\\ The
inequality (**) is clear for the canonical embeddings of the
projective space, namely the rank one case. Assume therefore that
$\rank G/K>1$. Observe that the embeddings $f_d$ factor through
the Veronese embeddings and the first canonical embedding, i.e.,
$f_d = \Ver \circ f_1$, where $\Ver: \C P^{N_1} \rightarrow \C
P^{N_d}$ is the Veronese embedding (see \cite[p.659]{NaTa} or
\cite[Section 3]{Take}). Then, the dimension of the normal space
of $f_d$ is greater than $\dfrac{N_1(N_1 + 1)}{2}$. Now recall
that any canonical embedding is full.  Thus, $m < N_1$ and
$$\dim_{\C} N_p > \dfrac{N_1(N_1 + 1)}{2} > \dfrac{m(m+1)}2 \,.
$$
 \hfill $\Box$

\begin{rem} By using the Weyl's formula \cite[p.189, Remark
2.3]{Take} to compute the dimension $N$ of the target projective
space $\C P^N$, it is possible to extend the above proof to an
arbitrary immersion of a homogeneous K\"ahler manifold $M$.
Namely, it is enough to check that the codimension of the
immersion $M \rightarrow \C P^N$ is bigger than $\dfrac{n(n +
1)}{2}$, $n=\dim_{\C}(M)$ to conclude that the normal holonomy
group is the full unitary group $U(N-n)$. This observation is
another motivation of Conjecture~\ref{ProBer}.
\end{rem}

\section{Complex parallel submanifolds.}\label{classification}

The goal of this section is to simplify the arguments in the
classical article \cite{NaTa}. Namely, we are going to avoid the
use of the eigenvalues of the curvature operator computed by
Calabi-Vesentini \cite{CaVe} and Borel \cite{Borel} which were
strongly used in \cite{NaTa}.

Here we are going to give a direct proof of the following theorem.
\begin{thm} \label{para0} Let $f_d : G/K \rightarrow \C P^{N_d}$ be
the $d$-th canonical embedding of an irreducible Hermitian symmetric
space $G/K$. Assume that the embedding has parallel
second fundamental form. Then, if $f_d$ is not the Veronese embedding, $f_d = f_1$, that is to say $f_d$
is the first canonical embedding.
\end{thm}

\it Proof. \rm As we remarked in the proof of
Theorem~\ref{alto},  the
embeddings $f_d$ can be described in terms of the Veronese embeddings and the first
canonical embedding. Namely,  $f_d = \Ver \circ f_1$, where $\Ver: \C
P^{N_1} \rightarrow \C P^{N_d}$ is the Veronese embedding. Then
the codimension of the embedding $f_d$ is greater than
$\frac{N_1(N_1 + 1)}{2}$ (and one has equality if any only if
$f_d$ is the Veronese embedding). Thus, if the codimension of
$f_d(G/K)$ is one then $f_d$ is the first canonical embedding of
$Gr_2^+(\R ^{n})$ i.e. a complex quadric. Thus, we can assume that
the codimension is greater than one and that $f_d$ is not the Veronese embedding.\\
Now recall that any canonical embedding is full. Let $n =
\dim_{\C}(G/K)$ be the complex dimension of $G/K$. Thus, $n < N_1$
and we get that the dimension of the first normal space is smaller or equal to
$\frac{N_1(N_1 + 1)}{2}$. Thus, if $d>1$ then $f_d$ cannot be
full since the first normal space is invariant by parallel
transport in the normal connection and thus agrees with the normal
space, by reduction of the codimension (i.e. by Theorem \ref{codi}).
This shows that $f_d = f_1$ and we are done. \hfill $\Box$\\

According to the above theorem, in order to get the embeddings of
Hermitian spaces with parallel second fundamental form we have to
study the first canonical embeddings only. The following theorem
gives a sharp description.

\begin{thm} \label{para}  Assume that the first canonical embedding of an
irreducible Hermitian symmetric space $M$ of higher rank has
parallel second fundamental form. Then $\rank(M)=2$.
\end{thm}

\it Proof. \rm  According to Theorem \ref{norhol}, if $M=G/K
\hookrightarrow \C P^N$ has parallel second fundamental form then
$K$ (or a quotient $S$) must act on the normal space $N_p$ as the
isotropy of an irreducible Hermitian symmetric space. From the
classification of irreducible Hermitian symmetric space we will
show that if $\rank(M) > 2$ then there is no irreducible Hermitian
symmetric space $H/S$ of dimension $\dim(N_p)$.
This will follow from a case by case analysis on the list in  the first column of
Table 2.\\
\begin{table} [h!]\ {\footnotesize
\vskip .5cm {}
\begin{tabular}{ | c | c | c | c | }\hline
\rule[-20pt]{0pt}{33pt}%
\parbox{3 cm}{\center{Hermitian symmetric space $G/K$}}  & \parbox{1.5cm}{\center{$\dim_\C G/K$}}
& \parbox{3.5cm}{\center {Codimension of its first canonical embedding}} &
\parbox{1.2cm}{\center {Rank of $G/K$}}\\ \hline\hline
\rule [-14pt]{0pt}{33pt}%
$\dfrac{E_7}{T^1 \cdot E_6}$ & $27$
& $28$ &  $3$ \\ \hline
 \rule [-14pt]{0pt}{33pt}%
$\dfrac{E_6}{T^1 \cdot Spin_{10}}$ & $16$  &
$10$ &  $2$ \\ \hline
 \rule [-14pt]{0pt}{33pt}%
 $\dfrac{Sp(n)}{U(n)}$ & $\dfrac{n(n+1)}2$  & $\left(\begin{array}{c} 2n\\n \end{array}\right)-\left(\begin{array}{c} 2n\\n-2  \end{array}\right)-1-\dfrac{n(n+1)}2$& $n$
  \\ \hline
 \rule [-14pt]{0pt}{33pt}%
 $\dfrac{SO(n+2)}{T^1 \cdot SO(n)}$ & $n$  & $1$ & $2$
  \\ \hline
 \rule [-14pt]{0pt}{33pt}%
 $\dfrac{SO(2n)}{U(n)}$ & $\dfrac{n(n-1)}2$  & $2^{n-1}-\dfrac{n(n-1)}2-1$& $[n/2]$ \\ \hline
 \rule [-14pt]{0pt}{33pt}%
$\dfrac{SU(a+b)}{S(U(a)\times U(b))}$ & $ab$  & $\left(\begin{array}{c} a+b\\ b \end{array}\right)-ab-1$ & min$(a,b)$ \\ \hline
\end{tabular}
} \caption{\footnotesize Hermitian symmetric spaces, their
dimensions, ranks and the codimension of their first canonical
embedding.}
\end{table} So, let us start with the rank $3$ exceptional
Hermitian symmetric space $\dfrac{E_7}{T^1 \cdot E_6}$. Notice
that the codimension of its first canonical embedding is $28$ (see
the third column of Table 2 constructed using \cite[p.
654]{NaTa}). Thus, a simple inspection on the second column of
Table 2 implies that there is no irreducible $28$-dimensional
Hermitian symmetric space whose isotropy is a quotient of $T^1
\cdot E_6$. Then, the first canonical embedding of
$\dfrac{E_7}{T^1 \cdot E_6}$ does not have parallel second
fundamental
form. \\
Going on on the list, $\dfrac{E_6}{T^1 \cdot Spin_{10}}$ is of rank 2, so we do not need consider it.\\
Let us consider further  $\dfrac{Sp(n)}{U(n)}$: the codimension of
its first canonical embedding is $h(n)=\left(\begin{array}{c} 2n\\n \end{array}\right)-\left(\begin{array}{c} 2n\\n-2  \end{array}\right)-1-\dfrac{n(n+1)}2$. There are two candidates for Hermitian
symmetric spaces whose isotropy is $U(n)$: $\dfrac{Sp(n)}{U(n)}$ and $\dfrac{SO(2n)}{U(n)}$.
However, some computations show that  dimension cannot be equal to $h$, for any $n$.\\
We can skip the case of $\dfrac{SO(n+2)}{T^1 \cdot SO(n)}$, since
its rank is two.\\ The Hermitian symmetric space
$\dfrac{SO(2n)}{U(n)}$ has first canonical embedding of
codimension $2^{n-1}-\dfrac{n(n-1)}2-1$. Again Hermitian symmetric
spaces whose isotropy is $U(n)$ are given by $\dfrac{Sp(n)}{U(n)}$
and $\dfrac{SO(2n)}{U(n)}$, but their dimensions cannot equal
$2^{n-1}-\dfrac{n(n-1)}2-1$ for any $n$. \\ Finally, for
$\dfrac{SU(a+b)}{S(U(a)\times U(b))}$ the codimension of its first
canonical embedding is $\left(\begin{array}{c} a+b\\ b
\end{array}\right)-ab-1$. Hermitian symmetric spaces whose
isotropy is a quotient of $S(U(a)\times U(b))$ are
$\dfrac{Sp(n)}{U(n)}$,  $\dfrac{SO(2n)}{U(n)}$ with $n=a$ or $n=b$
and $\dfrac{SU(a+b)}{S(U(a)\times U(b))}$, but none of them fits.
\hfill $\Box$\\

We are now going to show that the converse of Theorem \ref{para}
also holds. Namely,

\begin{thm} \label{para2}  The first canonical embedding of an
irreducible Hermitian symmetric space $M$ of rank two has parallel
second fundamental form. \end{thm}

\it Proof. \rm A first proof can be obtained by an explicit
construction of an extrinsic symmetry at each point of $M$. A
second proof was given in \cite[p.245]{Mok} by means of a
pinching theorem due to A. Ros \cite{Ros}. Another proof is in
Nakagawa-Takagi's paper \cite{NaTa}. \hfill $\Box$

\vspace{.5cm}

The above theorem is also a consequence of the following
conceptual argument. Indeed, one can see that the list of
submanifolds given by the images of the first canonical embedding
of an irreducible Hermitian symmetric space of rank two agrees
with the list of the unique complex orbits of the isotropy action
on the projective space $\P (T_{[K]}G/K)$ i.e. the second column
of Table 1. Thus, we just need to show the following proposition.

\begin{prop} \label{para1} The complex orbit of the (projectivized) isotropy
representation of an irreducible Hermitian symmetric space $G/K$
has parallel second fundamental form. \end{prop}

\it Proof. \rm
Let ${\g} = \k \oplus \p$ be a Cartan decomposition of the Lie
algebra $\g$ of $G$. Then, the isotropy action agrees with the restriction
to $\p\cong \C^{N+1}$ of the adjoint action of $K$.  \\
Let $M=K/K_0$ the unique complex orbit of the isotropy action on
the projective space $\C P^N$ and suppose that $M$ is an
irreducible Hermitian symmetric space. Let $\k=\k_0+\m$ be a
Cartan decomposition of $\k$. Then $M$ is isometrically embedded
in $\C P^N$ if it is endowed with a multiple of the opposite of
the Killing form on $\m$ and $\m\cong T_{[p]}M$. A similar
computation as in \cite[Lemma 4.1.5]{BCO} yields
$$[ \m, N_{[p]} M] \subseteq  T_{[p]}M\, . $$
This implies that $M$  has parallel second fundamental form, since
it is Hermitian symmetric (cf. \cite[Lemma 7.2.6]{BCO}). \hfill
$\Box$

\vspace{.5cm}

This completes the proof of Theorem~\ref{main}.

\begin{rem} Observe that Proposition \ref{para1} can be also proven
by using ideas in \cite[Chapter 6]{Mok}\ about characteristic
varieties. Furthermore, these orbits can also be described by using
the Jordan Algebra approach to Hermitian symmetric spaces, namely,
in terms of the so called \it tripotents \rm \cite{Roos} and
\cite[p. 579]{Kaup}.
\end{rem}

\begin{rem} Finally, we recall a problem proposed by A. Ros (Problem 5 in \cite[p. 272]{RosL}),
namely to characterize the symmetric submanifolds of $\C P^N$
without using metric notions (that is to say, changing
``holomorphic isometry'' with ``holomorphic transformation''). We
remark that Ros' problem could be related to the geometric
characterization of the so called \textit {Helwig spaces} see
\cite[p. 58, Problem II.4.5.]{Ber}.
\end{rem}

 {\small

 \noindent {\bf Acknowledgments.} We wish to thank Carlos Olmos for his many valuable suggestions.

\vspace{.5cm}

\noindent Dipartimento di Matematica, Politecnico di Torino,
\\ Corso Duca degli Abruzzi 24, 10129  Torino, Italy. \\ \tt
antonio.discala@polito.it \rm

\noindent Dipartimento di Matematica, Universit\`a di Torino, \\
via Carlo Alberto 10
10123 Torino, Italy \\
\tt sergio.console@unito.it
 \rm

}


\medskip


\begin{thebibliography}{12345}

\bibitem[AlDi]{AlDi}
{\sc Alekseevsky, D.V. and Di Scala, A.J.:} {\it The normal
holonomy group of K\"ahler submanifolds,} Proc. London Math. Soc.
(3) 89 (2004) 193-216.

\bibitem[BCO]{BCO}
{\sc Berndt, J., Console, S. and Olmos, C.:} {\it Submanifolds
and holonomy,} Research Notes in Mathematics 434, Chapman \&
Hall/CRC, 2003.

\bibitem[Ber]{Ber} {\sc Bertram, W.:}
{\it The geometry of Jordan and Lie structures}, Lecture Notes in
Mathematics 1754, Springer, 2000.

\bibitem[Borel]{Borel}
{\sc Borel, A.:} {\it On the Curvature Tensor of the Hermitian
Symmetric Manifolds,} Ann. of Math. (2) 71 (1960), no. 3,
508--521.

\bibitem[BoWe]{BoWe}
{\sc Borel, A. and Weil, A.:}  Repr\'esentations lin\'eaires et espaces homog\'enes k\"ahleriens
des groupes de Lie compacts, S\'eminaire Bourbaki, expos\'e no.
100 par J.-P. Serre,
(mai 1954).

\bibitem[Cal]{Cal}
{\sc Calabi, E.:} {\it Isometric imbeddings of complex manifolds,}
Ann. of Math. (2) 58 (1953), 1--23.

\bibitem[CaVe]{CaVe}
{\sc Calabi, E. and Vesentini, E.:} {\it On compact, Locally
symmetric Kaehler manifolds,} Ann. of Math. (2) 71 (1960), no. 3,
472--507.

\bibitem[Ce]{Ce}
{\sc Cecil, T.E.:} {\it Geometric applications of critical point
theory to submanifolds of complex projective space ,} Nagoya Math.
J. {\bf 55}, 5-31 (1974).


\bibitem[ChO1]{ChO1}
{\sc Chen, B-Y. and  Ogiue, K.:} {\it Some extrinsic results for
\Ka submanifolds ,} Tamkang J. Math. {\bf 4}  No. 2, 207-213
(1973).

\bibitem[DiOl]{DiOl}
{\sc Di Scala, A.J.  and Olmos, C.:} {\it Submanifolds with curvature
normals of constant length and the Gauss map,} J. Reine Angew.
Math. 574 (2004), 79--102.

\bibitem[Esch]{Esch}
{\sc Eschenburg, J.-H.:} {\it Lecture Notes on Symmetric Spaces,}
Preprint Augsburg 1998.


\bibitem[Fer]{Fer}
{\sc Ferus, D.:} {\it Symmetric submanifolds of Euclidean space,}
Math. Ann. 247 (1980), no. 1, 81--93.

\bibitem[GuSt]{GuSt}
{\sc Guillemin, V.  and Sternberg, S.:} {\it Symplectic techniques in
physics,} Cambridge University Press, Cambridge, 1984.

\bibitem[HeOl]{HeOl}
{\sc Heintze, E.  and  Olmos, C.:} {\it Normal holonomy groups and
$s$-representations,} Indiana Univ. Math. J. 41 (1992), no. 3,
869--874.

\bibitem[Kaup]{Kaup}
{\sc Kaup, W.:} {\it On Grassmannians associated with ${\rm JB}\sp
*$-triples,} Math. Z. 236 (2001), no. 3, 567--584.


\bibitem[Mok]{Mok}
{\sc Mok, N.:} {\it Metric Rigidity Theorems on Hermitian Locally
Symmetric Spaces,} Series in Pure Mathematics -Volume 6. World
Scientific, 1989.

\bibitem[Na]{Na}
{\sc Nakano, S.:} {\it On complex analytic vector bundles,}  J.
Math. Soc. Japan 7 (1955), 1--12.

\bibitem[NaTa]{NaTa}
{\sc Nakagawa, H.  and  Takagi, R.:} {\it On locally symmetric Kaehler
submanifolds in a complex projective space,} J. Math. Soc. Japan
{\bf 28} (1976), 638--667.

\bibitem[Ol1]{Olmos}
{\sc Olmos, C.:} {\it The normal holonomy group,} Proc. Amer.
Math. Soc. 110 (1990) 813-818.


\bibitem[Ol2]{Olmos2}
{\sc Olmos, C.:} {\it A geometric proof of Berger holonomy
theorem,} Ann. of Math. (2) 161 (2005), no. 1, 579--588.


\bibitem[Ol3]{Olmos3}
{\sc Olmos, C.:} {\it On the geometry of holonomy systems,}
Enseign. Math. (2) 51 (2005), no. 3-4, 335--349.

\bibitem[Ol4]{Olmos4}
{\sc Olmos, C.:} {\it Isoparametric submanifolds and their
homogeneous structures,} J. Differential Geom. 38 (1993), no. 2,
225--234.

\bibitem[Ol5]{Olmos5}
{\sc Olmos, C.:} {\it Homogeneous submanifolds of higher rank and
parallel mean curvature,} J. Differential Geom. 39 (1994), no. 3,
605--627.


\bibitem[Ros1]{Ros}
{\sc Ros, A.:} {\it A Characterization of Seven Compact Kaehler
Submanifolds by Holomorphic Pinching}, Ann. of Math. (2) 121
(1985), 377--382.

\bibitem[Ros2]{RosL}
{\sc Ros, A.:} {\it Kaehler submanifolds in the complex projective
space,} Lecture Notes in Math., 1209, Springer, Berlin, (1986),
259--274.




\bibitem [Roos] {Roos} {\sc Roos, G.:} {\it Jordan triple systems,}
 pp. 425-534, in \textit{J.~Faraut, S.~Kaneyuki, A.~Kor\'{a}nyi, Q.k.~Lu, G.~Roos,
Analysis and Geometry on Complex Homogeneous Domains}, Progress in
Mathematics, vol.\textbf{185}, Birkh\"{a}user, Boston, 2000.


\bibitem[Take]{Take}
{\sc Takeuchi, M.:} {\it Homogeneous K\"ahler submanifolds in
complex projective spaces,} Japan. J. Math. (N.S.) 4 (1978), no.
1, 171--219.


\end{thebibliography}
\end{document}